# A Straightforward Solution to Burnside's Problem
## S. Bachmuth

### 1. **Introduction**

The Burnside Problem for groups asks whether a finitely generated group, all of whose elements have bounded order, is finite. We present a straightforward proof showing that the 2-generator Burnside groups of prime power exponent are solvable and therefore finite. This proof is straightforward in that it does not rely on induced maps as in [2], but it is strongly dependent on Theorem B in the joint paper with H. A. Heilbronn and H. Y. Mochizuki [9]. Theorem B is reformulated here as Lemma 3(i) in Section 2. Our use of Lemma 3(i) is indispensable.

Throughout this paper we fix a prime power $q = p^e$ and unless specifically mentioned otherwise, all groups are 2-generator. At appropriate places, we may require $e = 1$ so that $q = p$ is prime; otherwise q may be any fixed prime power.

The only (published) positive results of finiteness of Burnside groups of prime power exponents are for exponents $q = 2, 3$ and $4$ ([10],[12]). Some authors, beginning with P. S. Novikoff and S. I. Adian {15}, (see also [1]), have claimed that groups of exponent k are infinite for k sufficiently large. Our result here, as in [2], is at odds with this claim.

Since this proof avoids the use of induced maps, Section 4 of [2] has been rewritten. Sections 2, 3 and 6 have been left unaltered apart from minor, mostly expository, changes and renumbering of items. The introduction and Section 5 have been rewritten. Sections 5 & 6 are not involved in the proof although Section 5 is strongly recommended. That section was intended to shed light on the Burnside groups, but it also hints at a simpler proof of the problem for primes. For example, the proof of Lemma 3(i) in [9[ for $q = p^e$ is an induction on e with $e = 1$ the start of the induction.

This paper has been written with the aim that familiarity or consultation with [2] is not necessary. Our belief is that this paper is self contained but for the necessity of reference [9], and to a lesser extent reference [3] for Theorem 1(ii). The proof of the crucial Lemma 3(i) in Section 2, (which is Theorem B and the Proposition on page 245 in [9]), occupies pages 239 -245 of [9]. The proof of Theorem 1(ii) is deduced from Theorem 1 of [3] which involves a matrix representation of the group of IA-automorphisms of a finitely generated free metabelian group. One need only restrict the representation to the group of inner automorphisms. In our restriction to two generators, this is unnecessary since it is also shown in [3] that the IA-automorphism group consists only of inner automorphisms. (A simple proof of this can be found in [7] (page 397, Corollary 3) and this result has been greatly generalized in [8].) However, Vaughan-Lee has used results from this paper to produce a self contained proof of Theorem 1(ii) - with no need for further references. (To the best of my knowledge this proof is unpublished. Although easy to explain to experts, Vaughan-Lee's proof will not be reproduced here since a proof for a general audience would needlessly lengthen this paper.) Otherwise, apart from reference [9], a general knowledge of groups and algebra should suffice for comprehension of this paper.

We begin the proof with a four part theorem which constructs a free group of matrices and three of its images. Details and proofs of all parts of Theorem 1 are in Sections 2 and 3.

Theorem 1

(i) We define a (rank 2) free group $F(\mathcal{R}[t, t^{-1}])$ consisting of 2x2 matrices having entries in the Laurent polynomial ring $\mathcal{R}[t, t^{-1}]$, where $\mathcal{R}$ itself is the Laurent polynomial ring with integer coefficients $\mathcal{R} = Z[x, x^{-1}, y, y^{-1}]$.

The letter F by itself will denote either the free group $F(\mathcal{R}[t, t^{-1}])$ or a generic free group. If it is necessary to distinguish which, the context will make it clear.

(ii) The (homomorphic) image $F(\mathcal{R})$ of $F(\mathcal{R}[t, t^{-1}])$, obtained by setting $t = 1$ in each element of $F(\mathcal{R}[t, t^{-1}])$, is a group of matrices with entries in $\mathcal{R}$ that is isomorphic to the free metabelian group $F/F''$ of rank 2.

(iii) For each positive integer n, there exists a quotient ring $S = S(n)$ of $\mathcal{R}$ such that if n is a prime-power, $F(S)$ is isomorphic to $F/F''F^n$, the free (rank 2) metabelian Burnside group of exponent n.

(iv) $F(S[t, t^{-1}])$ is solvable if (and only if) the integer n is a prime power.

We thus have the commutative square where the horizontal maps $\alpha$ and $\underline{\alpha}$ come from the ring homomorphism $\mathcal{R} \to S$ and the vertical maps are determined by sending t to the identity.

$$\begin{array}{ccc} F(\mathcal{R}[t, t^{-1}]) & \xrightarrow{\alpha} & F(S[t, t^{-1}]) \\ \downarrow & & \downarrow \\ F(\mathcal{R}) & \xrightarrow{\underline{\alpha}} & F(S) \end{array}$$

Note that for prime-power exponents, $\underline{\alpha}: F(\mathcal{R}) \to F(S)$ is the Burnside map on the free metabelian group. It is likely also true for arbitrary exponents, but our proof is valid only for prime-power exponents. (The Burnside map on the free group is $\beta: F(\mathcal{R}[t, t]) \to G$, where $G \cong F/F^q$.)

Although we use only Theorem 1(iv), as stated in this paper, the slightly more general statement formulated in the following corollary is worth recording.

Corollary: Suppose the surjective map $\gamma : F(\mathcal{R}[t, t^{-1}]) \to G$ is determined by the ring map $\mathcal{R} \to S$. Then G (the image of $\gamma$) is solvable.

Theorem 1(iv) is dependent on Theorem B in [9] and leads to a relatively easy solution of the Burnside problem in Section 4. We state this as Theorem 2.

Theorem 2: The Burnside group $F/F^q$ is an image of $F(S[t, t^{-1}])$ and hence solvable.

Since $F/F^q$ has only elements of finite order, it is a finite group as a consequence of Theorem 2. Not all solvable groups arising as images of $F(S[t, t^{-1}])$ are finite since some possess elements of infinite order. An example is $F(S[t, t^{-1}])$ itself, but whether finite or infinite, these groups have much in common with exponent q groups. All have at least one generator of order q and all have a commutator subgroup of exponent q. Section 5 will illustrate this in the study of $F(S[t, t^{-1}])$. The presence of such a large subgroup of exponent q in these groups lends a plausible explanation for their solvability.

The proofs of parts (i), (ii), (iii) of Theorem 1 are in Section 2, the proof of part (iv) is in Section 3. Section 4 contains the proof of Theorem 2. Section 5 examines properties of $F(S[t, t^{-1}])$. Section 6 considers the extension of Theorem 2 to finitely generated groups. Sections 2, 3, and 6 are similar, if not virtually identical, to the same sections in [2] in all essential details.

**Acknowledgment:** I am most grateful to George Bergman for discovering and helping correct several errors prior to submission of this manuscript. Most important, he uncovered an error in the main result of the previous manuscript [2] which is here explained. The following assumes knowledge of [2], but can and should be ignored if one has not read [2].

Theorem GB (Generalized Burnside) in Section 4 of [2] states:
Let $\gamma$ be a mapping of the free group $F(\mathcal{R}[t, t^{-1}])$ onto a group G such that $\gamma$ induces $\underline{\alpha}$. Then G (i.e., the image of $\gamma$) is solvable.

A counterexample to this statement is easily constructed by letting G be a direct product of a group which satisfies the hypothesis of the Theorem (e.g., a Burnside group) and a finite simple group. This direct product can be set up to satisfy the hypothesis, but would not be solvable. (The proof fails at the start of Proposition 2 in [2] since the fact that $M_1$ is in both $F(\mathcal{R})$ and $F(\mathcal{R}[t, t^{-1}])$ and hence implies that $\gamma(M_1) = \underline{\alpha}(M_1)$ is valid for only one of the coordinates in a direct product. Consequently an additional hypothesis to assure that Proposition 2 remains valid is required. Probably the easiest correction of Theorem GB is not to permit direct products, as follows:

Theorem GB: (Generalized Burnside): Let $\gamma$ be a mapping of the free group $F(\mathcal{R}[t, t^{-1}])$ onto a group G which is not a direct product. If $\gamma$ induces $\underline{\alpha}$, then G (i.e., the image of $\gamma$) is solvable.

## 2. Metabelian matrix groups

Let F be the free group of rank 2 generated by the 2x2 matrices $M_1$ and $M_2T$ with entries in the Laurent polynomial ring $Z[x, x^{-1}, y, y^{-1}, t, t^{-1}]$.

$$M_1 = \begin{vmatrix} 1 & 1-y \\ 0 & x \end{vmatrix} \qquad M_2T = \begin{vmatrix} yt & 0 \\ 1-xt & 1 \end{vmatrix}$$

That F is free can be seen by letting $x = 1$ and $y = t = -1$ and appealing to a well known theorem. cf. Sanov [16].

Notation: We set $\mathcal{R} = Z[x, x^{-1}, y, y^{-1}]$ and write $F = F(\mathcal{R}[t, t^{-1}])$. The units of $\mathcal{R}$ are the monomials $x^i y^j$ and $-x^i y^j$ where $i, j$ are integers. We call $x^i y^j$ the positive units.

Lemma 1: (i) The image of F upon setting $t = 1$ is isomorphic to the free metabelian group $F/F''$ generated by the matrices $M_1$ and $M_2$ with entries in $\mathcal{R} = Z[x, x^{-1}, y, y^{-1}]$.

We put $F(\mathcal{R}) = gp\langle M_1, M_2 \rangle$.

$$M_1 = \begin{vmatrix} 1 & 1-y \\ 0 & x \end{vmatrix} \qquad M_2 = \begin{vmatrix} y & 0 \\ 1-x & 1 \end{vmatrix}$$

(ii) Each element M of $F(\mathcal{R})$ has the form $M = uI + N$, where I is the identity matrix, u is a (positive) unit in $\mathcal{R}$, and the 2x2 matrix N has the form

$$N = \begin{vmatrix} \lambda_1(1-x) & \lambda_1(1-y) \\ \lambda_2(1-x) & \lambda_2(1-y) \end{vmatrix},$$

where $\lambda_1, \lambda_2$ in $\mathcal{R}$ satisfy $\lambda_1(1-x) + \lambda_2(1-y) = 1-u$. We will omit the I and write $M = u + N$.

Proof: (i) $F(\mathcal{R})$ is isomorphic to the group of inner automorphisms of the rank 2 free metabelian group as described in [3].

Note that [3] provides a description of the group of all IA- automorphisms of the free metabelian group of arbitrary finite rank, a group which contains the inner automorphisms. For rank 2, all IA-automorphisms of the free metabelian group are inner [3]. Note however that for ranks > 2 the IA-automorphism groups are more complicated. In rank 3 it is infinitely generated [5], but for ranks > 3, although complicated, the IA-automorphism groups become finitely generated [6]. The generalization of Lemma 1(i) for larger ranks presents no difficulty (see Section 6). For these larger ranks one need restrict the IA- automorphisms to the inner automorphisms to recover the free metabelian group of rank n.

(ii) A proof is contained in [4], Section 2, Proposition (ii). Because of the importance of this lemma, we shall reproduce the proof in an appendix to this paper.

Notation: Let $v = (1-x, 1-y)$ and write N as above in the form $N = [\lambda_i v]$, where $\lambda_i v$ denotes the $i^{th}$ row of N.

The lower central series for a group G is defined inductively by $G_1 = G$ and for $j > 1$ the $j^{th}$ term is $G_j = [G, G_{j-1}]$.

Lemma 2: If $M = u+N$ is in $F(\mathcal{R})$, where $N = [\lambda_i v]$, then

(i) $M^n = u^n + (1+u+u^2+\cdots+u^{n-1})N$ for any integer $n > 1$.

(ii) If M is in the commutator subgroup of $F(\mathcal{R})$, then $\lambda_1$ and $\lambda_2$ are in $\Sigma$, the augmentation ideal of $\mathcal{R}$. More generally, if M is in the $j^{th}$ term of the lower central series of $F(\mathcal{R})$, then $\lambda_1$ and $\lambda_2$ are in $\Sigma^{j-1}$.

Proof: Note that $N^2 = (1-u)N$ since $\lambda_1(1-x) + \lambda_2(1-y) = (1-u)$.

(i) From Lemma 1(ii) and the fact that $N^2 = (1-u)N$, we have
$M^2 = (u+N)^2 = u^2 + 2uN + (1-u)N = u^2 + (1+u)N$. Now proceed by induction.
$M^n = M^{n-1}M = (u^{n-1} + (1+u+\cdots+u^{n-2})N)(u+N) = u^n + (1+u+u^2+\cdots+u^{n-1})N$, using $N^2 = (1-u)N$.

(ii) Choose a basis for $G = F(\mathcal{R})$ modulo $G_{j+1}$ and, using induction on j, compute using this basis. A convenient basis (H. Neumann [14], chapter 3, section 6) is the basic commutators that are not in the second derived group. Notice that if M is in the commutator subgroup, then $u = 1$ and $\lambda_1(1-x) + \lambda_2(1-y) = 0$. Hence $\lambda_1$ and $\lambda_2$ are in the augmentation ideal of $\mathcal{R}$. In fact, observe that $[M_2, M_1] = I + [\lambda_i v]$, where $\lambda_1 = -(1-y)$ and $\lambda_2 = (1-x)$. This is the first case and the start of the induction. It is straightforward to check that if the basic commutator $[M_2, M_1, M_2, \ldots, M_2, M_1, \ldots, M_1] = I + [\lambda_i v]$, where $M_1$ and $M_2$ occur a total of $a+1$ and $b+1$ times respectively, then $\lambda_1 = -(1-y)^{b+1}(1-x)^a$ and $\lambda_2 = (1-y)^b(1-x)^{a+1}$. Furthermore, conjugation of a commutator by $M = u + N$ becomes multiplication (of the non identity term) by u. Thus, Lemma 2(ii) follows.

Before proceeding we need to introduce more notation.

$q = p^e$ as always denotes a prime power and we let $\phi(q) = p^e - p^{e-1}$, the Euler $\phi$-function. Let $\mathcal{I}(n)$ be the n-cyclotomic ideal, the ideal in $\mathcal{R}$ generated by all elements $1+u+u^2+\cdots+u^{n-1}$ where u is a positive unit in $\mathcal{R}$. We denote by $S$ the quotient ring $\mathcal{R}/\mathcal{I}(n)\Sigma$. In Lemma 3(i) & (iii) we require that n be a prime-power.

Lemma 3: (i) $\Sigma^{e\phi(q)} \subseteq \mathcal{I}(q)$ in $\mathcal{R}/\mathcal{I}(q)$. Equality holds if (and only if) $e = 1$, that is, $q = p^e$ is prime. Furthermore $\Sigma^{e\phi(q)-1} \not\subseteq \mathcal{I}(q)$.

If $e \geq 2$, $0 \leq j \leq e-1$ and $k \geq e(p^e-p^{e-1}) - j(p^{e-1}-p^{e-2})$, then $p^j\Sigma^k \subseteq \mathcal{I}(q)$.

(ii) The group generated by the matrices $M_1$ and $M_2$ over $\mathcal{R}/\Sigma^c$ is isomorphic to $F/F''F_c$ the free metabelian nilpotent group of class c.

(iii) The group $F(S)$ generated by the matrices $M_1$ and $M_2$ over $\mathcal{R}/\mathcal{I}(q)\Sigma$ is isomorphic to $F/F''F^q$, the free metabelian Burnside group of prime power exponent q.

Proof: (i) This is part of Theorem B in [9] and the Proposition on page 241 of [9]. Theorem B is the case $j = 0$.

Arbitrary finite rank is not addressed here. For now we are concerned only with rank $r = 2$, but these results from [9] are valid for ranks $r \leq p+1$. For rank $p+2$, Lemma 3(i) is no longer valid. The power of the augmentation ideal needed to fall into the cyclotomic ideal increases.

(ii) This follows from Lemmas 1 and 2(ii). For clearly a matrix of $F(\mathcal{R})$ which is in $F(\mathcal{R})_c$ has entries congruent to I modulo $\Sigma^c$ by Lemma 2(ii). Conversely, suppose $M = u + [\lambda_i v] \equiv I$ mod $\Sigma^c$. Then $u = 1$, $\lambda_i v$ are in $\Sigma^c$ and $\lambda_1(1-x) + \lambda_2(1-y) = 0$. Thus $\lambda_1 = \lambda(1-y)$, $\lambda_2 = -\lambda(1-x)$ and $\lambda$ is in $\Sigma^{c-2}$. From the proof of Lemma 2(ii), it is clear that M is in $F(\mathcal{R})_c$.

(iii) This follows from Lemmas 1 and 2(i). For again it is clear by Lemma 2(i) that a matrix in $F(\mathcal{R})'$ which is a $q^{th}$ power has entries congruent to I modulo $\mathcal{I}(q)\Sigma$. Conversely, suppose $M = u + [\lambda_i v] \equiv I$ modulo $\mathcal{I}(q)\Sigma$. Assume first that M is in the commutator subgroup. Then $u = 1$, $\lambda_i$ are in $\mathcal{I}(q)$, and $\lambda_1(1-x) + \lambda_2(1-y) = 0$. Thus the $\lambda_i$ are in the augmentation ideal and again $\lambda_1 = \lambda(1-y)$, $\lambda_2 = -\lambda(1-x)$ and both $\lambda(1-y)$ and $\lambda(1-x)$ are in $\mathcal{I}(q)$. By Theorem B of [9] or the Theorem of Dark and Newell [11], we conclude that $\lambda(1-y)$ and $\lambda(1-x)$ are in $\Sigma^{e\phi(q)}$ and thus M is in $F(\mathcal{R})^q$.

Any element M in F($\mathcal{R}$) can be written M = $M_1^i M_2^j C$, where C is in the commutator subgroup of F($\mathcal{R}$). Then M = $x^i y^j + [\lambda_i v]$ for suitable $\lambda_i$. By the hypothesis, i and j must each be a multiple of q. Hence, by Lemma 2(i), $M_1^i$ and $M_2^j$ are each congruent to the identity. So we may assume M is in the commutator subgroup. This completes the proof of Lemma 3.

3. **Power series representations of F ; Completion of Theorem 1**

In this section we will prove the crucial part iv of Theorem 1. Recall that $S = S(n) = \mathcal{R}/I(n)\Sigma$. We already know that F($S$) is solvable (i.e., metabelian); but if n = q is a prime power, we shall show that $F(S[t, t^{-1}])$ is also solvable.

We begin with the free group $F = F(\mathcal{R}[t, t^{-1}]) = gp\langle M_1, M_2T \rangle$ defined at the beginning of the previous section. For computational purposes, we shall describe $F(\mathcal{R}[t, t^{-1}])$ as a power series in $(t-1)^k$, k = 0,1,2,..., where the coefficients are 2x2 matrices over $\mathcal{R} = Z[x, x^{-1}, y, y^{-1}]$. Namely, let

$$T = \begin{vmatrix} t & 0 \\ 1-t & 1 \end{vmatrix} \quad \text{and} \quad U = \begin{vmatrix} 1 & 0 \\ -1 & 0 \end{vmatrix}$$

Then,
$M_2T = M_2(I + (t-1)U) = M_2 + (t-1)M_2U$
$(M_2T)^{-1} = (I + (t-1)U)^{-1}M_2^{-1} = (I - (t-1)U + (t-1)^2U^2 - (t-1)^3U^3 + \cdots)M_2^{-1}$
$\quad = M_2^{-1} - (t-1)UM_2^{-1} + (t-1)^2UM_2^{-1} - (t-1)^3UM_2^{-1} + \cdots$

Since $F(\mathcal{R}[t, t^{-1}]) = gp\langle M_1, M_2T \rangle$, we can state formally

Lemma 4: Any element f in $F(\mathcal{R}[t, t^{-1}])$ has a representation of the form
(*) $\qquad f = M_f + (t-1)A_1 + (t-1)^2A_2 + (t-1)^3A_3 + \cdots$
where $M_f \in F(\mathcal{R}) = gp\langle M_1, M_2 \rangle$ and the $A_i$ as well as $M_f$ are matrices over $\mathcal{R}$.

We know from Lemma 1(i) that F/F'' is the image of F via the map $f \to M_f$. The next few lemmas contain further consequences of this description of $F(\mathcal{R}[t, t^{-1}])$.

Notation: Recall that $\Sigma$ denotes the augmentation ideal of $\mathcal{R}$. If all the entries of a matrix A are in $\Sigma^i$ we say that A is in $\Sigma^i$.

Lemma 5: Suppose $f \in F$ as in (*) is in F' the commutator subgroup of F. Then all the $A_i$ (the coefficients of $(t-1)^i$) are in $\Sigma$.

Proof: In the free group $F(\mathcal{R}[t, t^{-1}])$, set $x = y = 1$ (i.e. factor by the ideal $\Sigma$ in $\mathcal{R}$). What remains is an infinite cyclic group (generated by the matrix we called T at the beginning of this section). Hence, the commutator subgroup of F is in the kernel of this map and the result follows.

Remarks:
(1) This proof readily generalizes when the rank is greater than two. (For rank r, setting $x_1 = ... = x_r = 1$ produces a free abelian group of rank r-1. See Section 6.)

(2) If f is in F", the second derived group of $F(\mathcal{R}[t, t^{-1}])$, then one can show that $A_1$, the coefficient of t-1, is in $\Sigma^3$. However the remaining coefficients are in $\Sigma^2$. For the next lemma, we will be content to merely assume that the coefficients of f in F" are in $\Sigma$ as asserted in Lemma 5. Using Lemma 5 as a starting point, we show that as we move down the derived series of $F(\mathcal{R}[t, t^{-1}])$, we double the power of the augmentation ideal in the entries of the coefficients of the $(t-1)^i$. The proof of Lemma 6(ii) follows an elegant argument of Vaughan-Lee.

Lemma 6: (i) g is in F" if and only if g has the form $g = I + (t-1)A_1 + (t-1)^2 A_2 + (t-1)^3 A_3 + \cdots$
(ii) If g is in $F^{(k)}$, the $k^{th}$ derived group of F, (k > 1), and $d = 2^{k-2}$, then
$g = I + (t-1)^d A_d + (t-1)^{d+1} A_{d+1} + \cdots$ where $A_d, A_{d+1}, A_{d+2},...$ are in $\Sigma^d$.

Proof: (i) is a restatement of Lemma 1(i).

(ii) The proof is by induction on k. The start of our induction is k = 2 which is covered by Lemma 5 and part (i) above. Suppose therefore that the result is true for $k \geq 2$ and let $d = 2^{k-2}$.
Let g, h be in $F^{(k)}$
$g = I + (t-1)^d A_d + (t-1)^{d+1} A_{d+1} + ...$, $h = I + (t-1)^d B_d + (t-1)^{d+1} B_{d+1} + ...$, and consider [g , h]. By induction we assume that $A_i$, $B_i$ are in $\Sigma^d$ for $i \geq d$. Since $g^{-1}$, $h^{-1}$ are also in $F^{(k)}$, we we can write $g^{-1} = I + (t-1)^d C_d + (t-1)^{d+1} C_{d+1} + ...$, $h^{-1} = I + (t-1)^d D_d + (t-1)^{d+1} D_{d+1} + ...$, and $C_i$, $D_i$ are in $\Sigma^d$ for $i \geq d$. With the summations over $i \geq d$, let $A = \Sigma (t-1)^i A_i$, $B = \Sigma (t-1)^i B_i$, $C = \Sigma (t-1)^i C_i$, $D = \Sigma (t-1)^i D_i$.
Then $gg^{-1} = (I+A)(I+C) = I+A+C+AC$, and hence $A+C+AC = 0$. Similarly $B+D+BD=0$. Hence,
$[g,h] = (I+A)(I+B)(I+C)(I+D) =$
$I+A+B+C+D+ AB + AC + AD + BC + BD + CD + ABC + ABD + ACD + BCD + ABCD = I + AB + AD + BC + CD + ABC + ABD + ACD + BCD + ABCD$.
So if we expand [g , h] as a power series

$$[g,h] = I + (t-1)E_1 + (t-1)^2 E_2 + (t-1)^3 E_3 + \cdots,$$

then $E_i = 0$ for $i < 2d$, and for $i \geq 2d$, $E_i$ is a linear combination of products of two or more matrices from the set $\{A_d, A_{d+1},\ldots, B_d, B_{d+1},\ldots, C_d, C_{d+1},\ldots, D_d, D_{d+1},\ldots\}$. By induction, all matrices in this set are in $\Sigma^d$, hence the $E_i$ are in $\Sigma^{2d}$ for $i \geq 2d$. This completes the proof.

Lemma 6(ii) is enough to show that the group $F(S[t, t^{-1}])$ is solvable. However, this implies a bound for the derived length of $F(S[t, t^{-1}])$ one larger than best possible since we assumed that an element of F" has all coefficients of $(t-1)^i$ in $\Sigma$, when it is easy to see that they lie in $\Sigma^2$. In order to give the best possible bound, we will use Lemma 7 in place of Lemma 6.

Lemma 7: If $g$ is in $F^{(k)}$, the $k^{th}$ derived group of F, ($k > 1$), and $d = 2^{k-2}$, then $g = I + (t-1)^d A_d + (t-1)^{d+1} A_{d+1} + \cdots$ where $A_d, A_{d+1}, \ldots$ are in $\Sigma^{2d}$.

Proof: Because of Lemma 6(ii), we only have to show that an element in F" has all coefficients in $\Sigma^2$. To show this, we use Lemma 5 as a starting point and proceed as in the proof of Lemma 6(ii). The only difference is that the constant term of an element in F' is no longer the identity matrix, but now has the form $u + [\lambda_i v]$, where $\lambda_1$ and $\lambda_2$ are in $\Sigma$. Otherwise the proof proceeds exactly as in Lemma 6(ii). We omit the details.

For the rest of this section we require that $q = p^e$ is a prime power in order to apply Theorem B of [9] which is needed for Theorem 1(iv). Theorem B links the augmentation and cyclotomic ideals of $\mathcal{R}$. It is summarized in Section 2 (Lemma 3(i)).

Theorem 1(iv): If $n = q = p^e$ is a prime power, $F(S[t, t^{-1}])$ is solvable. The derived length of $F(S[t, t^{-1}])$ is at most k where $2^{k-1} \geq e(p^e - p^{e-1}) + 1$.

Proof: We apply Lemma 7 to our previous calculations in $F(\mathcal{R}[t, t^{-1}])$ before mapping across to $F(S[t, t^{-1}])$. Choose k large enough so that $2^{k-1} \geq e\phi(q) + 1 = e(p^e - p^{e-1}) + 1$. Then, applying Lemma 3(i), the coefficients of the $(t-1)^i$ are the zero matrix in an expansion of an element of the kth derived group of $F(S[t, t^{-1}])$. Thus the elements in the $k^{th}$ derived group are the identity. This completes the proof of Theorem 1(iv).

Comments:
(1) Since $\Sigma^{e\phi(q)-1} \not\subset I(q)$, it is easy to see that the value of k in Theorem 1(iv) is best possible.

(2) Dark and Newell [11] have shown that the exact nilpotentcy class of $F/F''F^q$ is $e(p^e-p^{e-1})$ when F has rank 2 and $q = p^e$ is a prime power. Their result can be used to give an alternative method for calculating the solvability class of $F(S[t, t^{-1}])$ when $r = 2$. For our purposes, it appears more straightforward to use Theorem B in Bachmuth, Heilbronn, and Mochizuki [9], a method which remains viable for $r \leq p+1$.

We end this section with the following obvious Corollary which had been used in [2]. It may be worthwhile to state this explicitly as it may prove useful in further work.

Corollary: Suppose the surjective map $\gamma : F(\mathcal{R}[t, t^{-1}]) \to G$ is the result of sending $\mathcal{R}$ into $S$. Then G (the image of $\gamma$) is solvable.

Proof: We already observed that for the special case $G = F(S[t, t^{-1}])$, the Corollary is Theorem 1(iv). For this choice of G, the map $F(\mathcal{R}[t, t^{-1}]) \to F(S[t, t^{-1}])$ sends $\mathcal{R}$ to $S$ (while fixing t). For general G, we have as hypothesis that the elements of $\mathcal{R}$ are sent to $S$ in the map $\gamma$. Thus, as in Theorem 1(iv), Theorem B in [9] assures that the elements in a large enough power of the augmentation ideal are sent to the zero of $S$ in the map from $F(\mathcal{R}[t, t^{-1}])$ to G. Thus, in the expansion of the derived series in G, the elements in a high enough solvability class become the identity in G.

Remarks:
1) In the proof of Theorem 1(iv), we are in effect leaving intact the procedure of Vaughan-Lee (as was done in [2]) which showed that $F(S[t, t^{-1}])$ is solvable. We use the hypothesis $\mathcal{R} \to S$ to show that the image of the map $F(\mathcal{R}[t, t^{-1}]) \to G$ is solvable in place of the original map $F(\mathcal{R}[t, t^{-1}]) \to F(S[t, t^{-1}])$, a map defined by $\mathcal{R} \to S$.

2) The image of t determines whether G is finite or infinite. Moreover, it plays a role in the determination of the solvability class. If t is left fixed, then $G = F(S[t, t^{-1}])$ is an infinite solvable group with the largest possible solvability class for the given q. Enlarging the ideal will generally decrease the solvability class. In the extreme case of sending t to 1, the ideal becomes as large as possible, and G becomes $F(S)$ which has the smallest possible solvability class, (i.e., metabelian or even abelian when $q = 2$).

3) Theorem 1(iv) solves a problem in groups (establishing solvability) via calculations in a ring. This procedure was prominent in Wilhelm Magnus' classic paper [13]. There and in further papers by various authors since, the terms in a higher commutator series are linked to a 'degree function' in a ring - in our case to the powers of the augmentation ideal. Theorem 1(iv) is achieved by showing that as one proceeds down the derived series of $F(\mathcal{R}[t, t^{-1}])$, the powers of the augmentation ideal of $\mathcal{R}$ increase. Theorem B in the paper with Heilbronn and Mochizuki [9] ensures that they eventually fall into the q-cyclotomic ideal of $\mathcal{R}$, that is, the zero of $S$. Thus, upon the transfer of $F(\mathcal{R}[t, t^{-1}]) \to G$, we are able to conclude that G is solvable. This is the type of argument employed by Magnus [13] in his solution of a problem posed by Hopf.

4. **Solution of the Burnside Problem**

The Burnside group (of exponent q) is the largest group such that $f^q = 1$ for all f in F. The Burnside Problem asks whether such a finitely generated group is finite. Our solution is for two generator groups of prime-power exponents. Section 6 discusses the finitely generated case. We begin with the following lemma.

Lemma 8: Suppose f lies in the second derived group of the free group $F = F(\mathcal{R}[t,t^{-1}])$ and has image g in $G = F(S[t,t^{-1}])$. Then the order of g is a divisor of q.

Proof: We mention that g has order $q = p^e$ unless g is a $p^i$ power of a commutator, $i \le e$, in which case g has order $p^{e-i}$. In what follows, we will repeatedly use the fact that the ring $S$ has characteristic q. Represent f as in Section 3, and note that since f is in the second derived group, $M_f$ is the identity. Thus, $f = I + N$, $N = (t-1)A_1 + (t-1)^2 A_2 + (t-1)^3 A_3 + \cdots$, where the $A_i$ are matrices over $\mathcal{R}$, and the $A_i$ are in $\Sigma$, the augmentation ideal (Lemma 5 in Section 3). (Although not necessary, we will simplify the computations by assuming that the $A_i$ are in $\Sigma^2$ since f is in the second derived group.) Suppose that $e = 1$, that is q is a prime p. Then $g^p = (I + N)^p = I + N^p$, where $N^p = ((t-1)A_1 + (t-1)^2 A_2 + (t-1)^3 A_3 + \cdots)^p = ((t-1)^p B_1 + (t-1)^{p+1} B_2 + \cdots$. Each of the $B_i$ have entries in the 2p power of the augmentation ideal and since $2p \ge p-1$, Lemma 3(i) implies that each $B_i$ is the zero matrix, whence $g^p = 1$. Assume next that $e = 2$, ($q = p^2$). Then omitting (unnecessary) units from the coefficients, $g^q = (I + N)^q = I + pN^p + N^q$. In a calculation similar to the prior ones, the matrix in the $(t-1)^j$ term of $N^q$ has entries in the 2q power of the augmentation ideal and since $2q \ge 2(p^2 - p)$, $N^q$ is zero. Similarly, the entries of the matrices in $pN^p$ all lie in the $2p^2$ power of the augmentation ideal. Since $2p^2 \ge 2(p^2-p) - (p-1)$, Lemma 3(i) implies that $pN^p$ is also zero. Thus $g^q = 1$. For general e, with the help of the binomial theorem, one proceeds similarly. For the prime power factors $p^j$ of q, we need to know the values of k for which $p^j \Sigma^k \subseteq \mathcal{I}(q)$. Lemma 3(i) contains the needed information, namely:

If $e \ge 2$, $0 \le j \le e-1$ and $k \ge e(p^e - p^{e-1}) - j(p^{e-1} - p^{e-2})$, then $p^j \Sigma^k \subseteq \mathcal{I}(q)$. We leave it for the reader to fill in the details for $e \ge 3$.

Remark: Section 5 has a different proof than the one in Lemma 8, which for simplicity is limited to primes. We mention in particular that the proof of the Burnside Problem for primes is seriously simplified in the all important Lemma 3(i).

**T**heorem 2 (Burnside's Problem):  $F/F^q$ is solvable (and hence finite) for any prime power q.

Proof:  Begin by mapping the free group $F = F(\mathcal{R}[t,t^{-1}])$ to $G = F(S[t,t^{-1}])$. This is a surjective map from F to the solvable group G. We shall examine the kernel of this map and conclude that the elements in the kernel lie in the normal subgroup of F generated by its qth powers.

The generators of F, $M_1$ and $M_2T$, are defined in Section 2. If we set $t = 1$ in $M_2T$ we get $M_2$ and if we set $x = y = 1$ in $M_2T$ we get T. Neither $M_2$ nor T are elements of the free group F, but when convenient we can consider the generator $M_2T$ as the product of the matrices $M_2$ and T. In this way it is clear what is meant by the exponential sum of T for a word in F. We need only set $x = y = 1$ in a word in F and observe the exponent of T that ensues. We next remind the reader that the first part of Section 2 defines the ring $S$ which endows the group G with numerous elements of exponent q. In particular if we eliminate T in F by setting $t = 1$, the map $F = F(\mathcal{R}[t,t^{-1}])$ to $G = F(S[t,t^{-1}])$ becomes the map $F(\mathcal{R})$ to $F(S)$. By Lemma 3(iii), the kernel of the map $F(\mathcal{R})$ to $F(S)$ is the normal subgroup generated by the qth powers of $F(\mathcal{R})$. This includes the qth powers of commutators of $F(\mathcal{R})$ which is already explicit in Lemma 2(i) (for commutators $u = 1$). Lemma 2(i) shows that an element of order q corresponds to a cyclotomic polynomial in $\mathcal{R}$. (This was the reason for choosing the cyclotomic ideal in $\mathcal{R}$ to define $S$.)  Finally we caution the reader that we use the same symbols x, y, t for the units of $\mathcal{R}$ or the units of its quotient ring $S$. The context should be clear as to when we are in F or in G. With these preliminaries dispensed, we examine the kernel of the map from F to G by examining the image of all words in F.

The generator $M_1$ of F has order q in G. Thus $M_1^q$ is in the kernel of our map, and powers of $M_1$ that are not multiples of q are not in the kernel. Consider next the second generator $M_2T$ of F. Since T is left fixed by the map from F to G, send $(M_2T)^k$ into G and then set $x = y = 1$. This leaves $T^k$ which has infinite order in G when $k \neq 0$. Hence these elements are not in the kernel for $k \neq 0$. More generally, let f be any element of F not in the normal closure of the subgroup generated by $M_1^q$. Since f is a product of the generators $M_1$ and $M_2T$, we again set $x = y = 1$ in f and again we are left with $T^k$ for some integer k. If $k \neq 0$, the same argument which shows that $(M_2T)^k$ is not in the kernel also shows that f is not in the kernel. We are thus left with the case $k = 0$. That is, the exponential sum of T is zero and therefore the exponential sum of $M_2T$ is also zero.  In this case f might be in the commutator subgroup. Since this is dealt with separately, assume that the exponential sum of T is 0, but f is not in the commutator subgroup. Since f is a product of the generators $M_1$ and $M_2T$, this implies that f is a product of transforms of powers of $M_1$.

Put f in the form $uI + V$, where I is the identity matrix, u is a unit $x^i y^j t^k$ in $\mathcal{R}[t,t^{-1}]$ and V a matrix over $\mathcal{R}[t,t^{-1}]$. (To put f in this form, write $M_1$ and $M_2$ as in Section 2, Lemma 1(ii) and write the matrix T (described in Section 3) in the form $tI + W$, where W is the matrix with 0 entries in the first row and t-1, -(t-1) in the second row. f is in the kernel if and only if u = 1 and V is the zero matrix in G. That is, $f = I + V$, where all entries of V are in the cyclotomic ideal of $\mathcal{R}$. The exponent of t in u is the exponent of T in f. If this exponent is non-zero, (i.e., $k \neq 0$), then f is not in the kernel. Thus for f to be in the kernel, k = 0 and then necessarily j = 0. We have earlier observed that if $u = x^i$, $i \neq 0$, then f is in the kernel if i = q since $x^q = 1$ in $S$. This happens if $f = M_1^q$ or if f is a transform of $M_1^q$. In the more general case when the exponential sum of T is 0, we observed that for f to be in the kernel, the exponential sum of $M_1$ must be a multiple of q. In this case u = 1 and we need check that V is the zero matrix. This is a deterministic procedure and we could leave it at that, i.e., f is in the kernel if and only if V = 0. However, it is easy to see, with the help of Lemma 3(i), that V will always be zero if u = 1. Namely, rewrite f as as product of transforms of $M_1^i$, where the sum of the i are congruent to 0 modulo q. Upon calculating V, similarly as done in Section 2, the entries have the same power of the augmentation ideal as when calculating the product of the $M_1^i$ and are thus in the cyclotomic ideal. Another way to view this is to set t = 0 which sends f to the free metabelian group $F(\mathcal{R})$. This enables one to use the results in Section 2 to conclude that this image of f is in the normal closure of the qth powers. If we pull back to f in $F(\mathcal{R}[t,t^{-1}])$, the conjugates of $M_1^i$ by T will not have any material effect on the power of the augmentation in the entries of V. Thus V = 0 and f is in the normal closure of the qth powers.

We are left with the case i = j = k = 0, that is, the element f of F has exponential sum 0 in both generators of F. In particular, f is in the commutator subgroup of F.

To investigate this case where the exponential sum of both generators is zero, it is convenient to represent f in F as in Section 3. This enables one to use methods and results from that section without the need for explanations. Thus f in F is represented by a series in powers of (t-1), $f = M_f + (t-1)A_1 + (t-1)^2 A_2 + (t-1)^3 A_3 + \cdots$, where $M_f \in F(\mathcal{R}) = gp\langle M_1, M_2 \rangle$ and the $A_i$ as well as $M_f$ are matrices over $\mathcal{R}$. For f to be in the kernel of the map from F to G we must have $M_f = 1$ and the entries in each $A_i$ the zero of $S$.

For f in the commutator subgroup of F, the entries of the $A_i$ are in the augmentation ideal of $\mathcal{R}$ (Lemma 5 in Section 3). By Lemma 8 we know that the qth power of elements in the second derived group are in the kernel. (It is also true for all qth powers in the derived group and a proof is given for prime exponents in Section 5.) The question considered now is if there are other elements in the kernel that are not necessarily qth powers.  Section 3, combined with Lemma 3(i) of Section 2, gives  a procedure for determining if an element in the commutator subgroup is in the kernel.  There are two possibilities that $M_f = 1$.  The first is that $M_f$ is a qth power in F(R) as considered earlier in this section.  If $M_f = M_1{}^q$, then f does not contain T and f is in the kernel.  If $M_f$ is a qth power that contains the generator $M_2 T$, then f contains the variable t and thus has infinite order.  The second case where $M_f = 1$ is when f is in the second derived group of F.  We now consider this possibility in the search for kernel elements.

Consider two elements of F having exponent q and thus both in the group generated by qth powers.  We may regard their commutator as a product of transforms of these elements.  Thus the commutator is in the normal subgroup generated by qth powers and, as observed, its entries are in the augmentation ideal.  In iterating this procedure, we remain in this normal subgroup and, as detailed in Section 3, the power of the augmentation ideal in the matrix entries increase.  Ultimately, for f in a sufficiently large (depending on q) term of the derived series of F, the entries of the $A_i$ achieve a large enough power of the augmentation ideal to land in the cyclotomic ideal by Lemma 3(i).  Furthermore,  Lemma 3(i) informs us when this power is large enough and therefore when each $A_i$ is the zero matrix.  If the entries of an element f are not in a sufficiently high power of the augmentation ideal for Lemma 3(i) to apply, then f is not in the kernel.  Our inductive procedure asserts that these kernel elements are in the normal closure of the group generated by qth powers.  (The order of such an element will be a proper divisor of q if the commutator itself is a power that is a proper divisor of q.  In such cases Lemma 3(i) will still determine if the entries are in the cyclotomic ideal - as was exhibited in the proof of Lemma 8.)

We have examined all possibilities and can summarize the content of the kernel K of F → G as follows:  Consider $F/K \cong G$.  If f in F does not contain T or has exponential sum zero in T, then  $f^q$ is in K.  If the exponential sum of T in f is not zero, then for any $j \neq 0$, $f^j$ is not in K.  Other elements in the kernel, which are not obviously qth  powers, lie in a large term of the derived series - large enough for their entries to be zero as determined by Lemma 3(i).  They are a linear combination of cyclotomic polynomials.  Although such expressions of cyclotomics are readily identifiable in $F(\mathcal{R})$, we are not able to readily identify them in $F(\mathcal{R}[t,t^{-1}])$.  However, being group elements, they correspond to products and/or transforms of exponent q elements. These are relations obtained from the qth powers.  Upon the conclusion of Theorem 2, we will compare these results with prior knowledge of Burnside Groups.

Thus the kernel K of F → G contains elements of order q (corresponding to generators of the cyclotomic ideal) or group elements corresponding to elements in the ideal generated by the cyclotomic polynomials.  In particular, K is contained in the normal closure of the qth powers of F.  It follows that the Burnside group of prime power exponent is an image of the solvable group G and hence is a solvable and therefore a finite group.  This completes the proof of Theorem 2.

We next state a few results derivable from this Theorem which can be compared with known results: Let u and v denote units in $\mathcal{R}/\mathcal{I}(q)$. By convention the abelian groups are solvable groups of class 1.

q = 2: The proof tells us that 1-u is a cyclotomic polynomial - which it obviously is. Thus Burnside groups of exponent 2 are abelian.

q = 3: Here $(1-u)^2$ is a cyclotomic polynomial, and (1-u)(1-v) is in the cyclotomic ideal. Thus Burnside groups of exponent 3 are metabelian.

q = 4: Here $2(2^2-2^1) = 4$, $2(2^2-2^1) - ((2^1-1) = 3$. Thus $(1-u)^4$, $(1-u)^i (1-v)^j$, i+j = 4, and $2(1-u)^3$, $2(1-u)^i (1-v)^j$, i+j = 3 are in the cyclotomic ideal. Hence Burnside groups of exponent 4 are solvable of class 3.

q = 5: Here $(1-u)^4$ is a cyclotomic polynomial, while $(1-u)^3(1-v)$, $(1-u)^2 (1-v)^2$, and $(1-u)(1-v)^3$ are in the cyclotomic ideal. Thus Burnside groups of exponent 5 are solvable of class 3.

q = 7: Here $(1-u)^6$ is a cyclotomic polynomial, while $(1-u)^5(1-v)$, $(1-u)^4 (1-v)^2$, $(1-u)^3 (1-v)^3$, etc. are in the cyclotomic ideal. Thus Burnside groups of exponent 7 are solvable of class 4.

At q = 5 we have passed the point where we are able to compare results with past knowledge. We end in this vein with an example of a "large" prime power.

$q = 5^3$: Here $3(5^3-5^2) = 300$, $3((5^3-5^2) - j((5^2-5) = 260$ if j = 2 and 280 if j = 1. Thus $(1-u)^{300}$, $25(1-u)^{260}$, $5(1-u)^{280}$, $(1-u)^i(1-v)^j$ where i + j = 300, $25(1-u)^i(1-v)^j$, i + j = 260 and $5(1-u)^i(1-v)^j$, i + j = 280. are all in the cyclotomic ideal. As a consequence Burnside groups of exponent 125 are solvable of class 9.

Expressions such as $(1-u)^{139}(1-v)^{161}$ are in the cyclotomic ideal and are used here to establish solvability of the Burnside groups. The proof of such results for any q can be found in Theorem B and the Proposition following Theorem B in [9]. The following remark is an interesting idea to explore as a guide in understanding the relationship of elements in F to the cyclotomic ideal:

Remark: Contrast the results and proofs of Theorem 2 with the situation when we eliminate T in F by setting t = 1. This latter procedure produces the free metabelian group $F(\mathcal{R})$ which was studied in Section 2. By sending $\mathcal{R}$ to $\mathcal{S}$ in $F(\mathcal{R})$ we get $F(\mathcal{S})$, the free metabelian Burnside group. We therefore have a situation similar to above, but without the t variable present. This enables one to compare the results of the map F → G with the map $F(\mathcal{R}) \to F(\mathcal{S})$. In both situations we use the same ideal in $\mathcal{R}$ generated by the cyclotomic polynomials. These polynomials correspond to qth powers in F or $F(\mathcal{R})$ respectively. This allows one to interpret an element in the ideal with the corresponding element in $F(\mathcal{R})$ as a help in the understanding of the Burnside group. In $F(\mathcal{R})$, using Lemma 1(ii) in Section 2, it is easy to ascertain that conjugation of a commutator by an element $M = u + N$ becomes multiplication (of the non identity term) by u. Furthermore, multiplication of transforms of qth powers correspond to addition of polynomials. In this manner we can associate elements in the ideal with elements in the kernel. Notice however that the kernel of $F(\mathcal{R}) \to F(S)$ is the normal subgroup generated by all qth powers, whereas the kernel of $F \to G$ is a proper subgroup of the normal subgroup generated by the qth powers. This is due to the presence of T. Leaving T unchanged sets G as an intermediary in the construction of the Burnside group. The similarities in the behavior of the maps $F \to G$ and $F(\mathcal{R}) \to F(S)$ merit further study.

As a final observation, we note that a consequence of the Burnside group being an image of G is that the Burnside group of prime exponent $q = p$ requires at least the relation $t^p = 1$ ( equivalently $(t-1)^p = 0$) added to the definition of G. We state this formally as a corollary. It is doubtless true for arbitrary q, but we restrict our verification to primes.

Corollary: The Burnside group of prime exponent p is an image of the quotient group $G = F(S[t, t^{-1}])$ and the additional relation $(t-1)^p = 0$ added to the ring $S[t, t^{-1}]$.

Proof: Let $f = M_2T$, the second generator of F. Write $f = M_2T = M_2 + (t-1)M_2U$. Notice that the matrix U satisfies $U^2 = U$. Map $f \to g$ via the map $F \to G$. Then $g = M_2 + (t-1)M_2U)$, where $M_2$ and U are now matrices over S. (The entries of U are no longer integers but integers mod p, and $M_2^p = 1$.) We have $g^p = (M_2 + (t-1)M_2 U)^p = 1 + (t-1)^p(M_2U)^p$. Using $U^p = U$, set $x = y = 1$ in $g^p$ to conclude that $g^p = 1$ if and only if $(t-1)^p = 0$.

We leave as a conjecture the question as to whether this defines the Burnside group of exponent p. That is, the group $G = F(S[t, t^{-1}])$ and the additional relation $(t-1)^p = 0$ added to the ring $S[t, t^{-1}]$.

Conjecture: The Burnside group of prime exponent q is the group $G = F(S[t, t^{-1}])$, $t^q = 1$.

Note: Adding the relation $(t-1)^q = 0$ amounts to enlarging the q-cyclotomic ideal in the definition of the ring S to include all cyclotomic polynomials $1+u+u^2+\cdots+u^{q-1}$ where u is a positive unit in the enlarged ring $\mathcal{R}[t, t^{-1}]$.

Prior to reading Section 6 for the more general problem of finitely generated groups, we remind the reader that Lemma 3(i) in Section 2 is valid up to $p + 1$ generator groups, where $q = p^e$, and is false for groups requiring more than p+1 generators. However, if $e = 1$, that is for q a prime, Lemma 3(i) is valid for any finite number of generators. Thus for primes, Section 6 may display a path for a description of the Burnside group having an arbitrary number of generators without the need to extend Lemma 3(i).

## 5. The group $G = G(q) = F(S[t, t^{-1}])$.

We saw that the Burnside group of prime power exponent is an image of $G = G(q) = F(S[t, t^{-1}])$. Further consideration of G and its homomorphic images may shed a better understanding of the Burnside groups.

The group $G = F(S[t,t^{-1}])$ depends on the ring $S = S(n)$ and we can write $G = G(n)$ if we wish to emphasize this dependence on n. As we have seen, when n is a prime power $q = p^e$, $G = G(q)$ is solvable and the results in [9] enable one to determine the solvability class. Obviously this class is determined by the ideal in $\mathcal{R}$ that is used to construct $S$, an ideal which depends on q. Many technical questions are easily answered. For example, it is trivial to see that the solvability class of G(q) tends to infinity with increasing q. However, it appears puzzling why in the first instance G(q) should be solvable for any q. What properties does G(q) possess that might predict solvability? This question motivates what follows and perhaps sheds some light on the Burnside Problem. To simplify matters in this section, we will assume q is a prime.

Let $S = S(q)$, where $q = p$ is a prime. The generators of G are the images of $M_1$ and $M_2T$ under $\alpha$. To simplify notation we will use the same letters for the elements of G as for $F(\mathcal{R}[t,t^{-1}])$ but with a line through them. Hence the lined letters are matrices over $S$. Thus $G = F(S[t,t^{-1}]) = $ gp $\langle \alpha(M_1), \alpha(M_2T) \rangle =$ gp $\langle \overline{M_1}, \overline{M_2}T \rangle$. T is left unlined since $\alpha$ leaves T unchanged.

In Section 4 we observed that $\overline{M_1}$ has order q and $\overline{M_2}T$ has infinite order. If $\overline{W}$ is an element of G, i.e., a product of $\overline{M_1}$, $\overline{M_2}T$ and their inverses, such that $\overline{W}$ has a nonzero exponential sum in T, then we saw that $\overline{W}$ has infinite order. Next consider the elements with zero exponential sum in T. An example is a conjugate of $\overline{M_1}$ which of course has order q. However, there is a better method to show this, a method which ultimately relies upon Theorem B in [9] (i,e., Lemma 3(i)). That proof shows that any element which has exponential sum zero in T has order q. We state this as:

Theorem 3: Assume that q is a prime. Let $\overline{W}$ be an element of $G = F(S[t,t^{-1}])$.
a) If $\overline{W}$ has exponential sum zero in T, then $\overline{W}^q = 1$.
b) If $\overline{W}$ has a nonzero exponential sum in T, then $\overline{W}$ has infinite order

Proof: If $\overline{W}$ has exponential sum zero in T, then we claim that $\overline{W}^q = 1$. $\overline{W}$ is a product of $\overline{M_1}$, $\overline{M_2}T$ and their inverses. Let W be a preimage of $\overline{W}$ in $F(\mathcal{R}[t,t^{-1}])$. Write the $M_i$ in the form of Lemma 1(ii) of Section 2 and write T in the form $tI + X$, where tI is the scalar matrix with t in the diagonal. As a product of these matrices W has the form $W = uI + V$, where $u = x^i y^j t^k$ and the entries of V are in the augmentation ideal of $\mathcal{R}$. By our assumption, k = 0 (the only place we need the hypothesis that the exponential sum of T is zero). Now map W to $\mathbb{W}$, i.e., send $\mathcal{R}$ to $S$. Since $S$ has prime characteristic q, $(uI + \mathbb{V})^q = u^q I + \mathbb{V}^q$, where $u^q = 1$ and $\mathbb{V}^q$ is the zero matrix. Namely, the entries of $\mathbb{V}^q$ are in the qth power of the augmentation ideal. By Lemma 3(i), the q-1 power of the augmentation ideal lies in $\mathcal{I}(q)$ and hence the qth power of $\mathbb{V}$ is in $\mathcal{R}/\mathcal{I}(q)\Sigma$ which is the zero of $S$.

We next combine Theorem 1 and Theorem 3. Recall that setting t = 1 in a matrix of G sends that matrix into one whose order is a divisor of q. (Factoring by the ideal t - 1 in $S[t,t^{-1}]$ is factoring by the second derived group in G.) Thus as a corollary of Theorems 1 and 3, we have the

Proposition: Let g be an element in $G = G(q) = F(S[t,t^{-1}])$. Then either $g^q = 1$ or g has infinite order. If g has infinite order, $g^q \cong 1$ modulo G''.

The ring $S = S(q) = \mathcal{R}/\mathcal{I}(n)\Sigma$ was chosen as precisely the ring which sends the free metabelian group $F(\mathcal{R})$ into $F(S)$, the Burnside metabelian group of exponent q. Using the same ring map to change $F(\mathcal{R}[t,t^{-1}])$ into $F(S[t,t^{-1}])$ one might expect to find a substantial number of elements of exponent q in $F(S[t,t^{-1}])$. Theorem 3 describes the ubiquity of these elements. This leads to a possible group theoretic explanation to the question of why $F(S[t,t^{-1}])$ and its images are solvable. Namely, enough elements of exponent q lead to (higher) commutator identities as in $F(S[t,t^{-1}])$. For small q such identities were known in groups of exponent q. The ideal in $\mathcal{R}$ used to construct $S$ was chosen to produce elements of exponent q. Then with the help of Theorem B in [9], commutator identities, seemingly unexpectedly, appear in these groups.

## 6. Arbitrary Generators

We examine the main ingredients necessary for generalization to k generators.
$\mathcal{R} = \mathcal{R}(k) = Z[x_1, x_1^{-1}, \ldots, x_k, x_k^{-1}]$ is the k generator Laurent polynomial ring over the integers.
$\mathcal{R}[t, t^{-1}] = Z[x_1, x_1^{-1}, \ldots, x_k, x_k^{-1}, t, t^{-1}]$ is the k+1 generator ring.
$\mathcal{I}(q)$ is the q-cyclotomic ideal in $\mathcal{R}$, the ideal generated by all q-cyclotomic elements $1 + u + u^2 + \cdots + u^{q-1}$, where u is a positive unit in $\mathcal{R}$ and q is the prime power $q = p^e$.

$S$ is the quotient ring $\mathcal{R}/\mathcal{I}(q)\Sigma$, where $\Sigma$ is the augmentation ideal of $\mathcal{R}$.
We now describe four k x k matrix groups: $F(\mathcal{R})$, $F(S)$, $F(\mathcal{R}[t, t^{-1}])$ and $F(S[t, t^{-1}])$.

Let v be the k-tuple, $v = (1-x_1, \ldots, 1-x_k)$ and N the k x k matrix $N = [\lambda_i v]$, where $\lambda_i v$ denotes the ith row of N, $i = 1,\ldots,k$. $F(\mathcal{R})$ is the set of all matrices $uI + N$, where I is the identity matrix, u is a (positive) unit in $\mathcal{R}$, and the $\lambda_1, \ldots, \lambda_k$ in $\mathcal{R}$ satisfy $\lambda_1(1-x_1) + \ldots + \lambda_k(1-x_k) = 1-u$. For a matrix M in $F(\mathcal{R})$ we often omit the I and write $M = u + N$. Alternatively, one may define $F(\mathcal{R})$ by the generators: $F(\mathcal{R}) = gp\langle M_1, \ldots, M_k \rangle$, where for $j=1,\ldots,k$ $M_j = x_j I + [\lambda_i v]$, $\lambda_i = 0$ unless $i = j$ in which case $\lambda_j = 1$.

Replacing $\mathcal{R}$ by $S$ defines $F(S)$ as well as $F(S[t, t^{-1}])$ as soon as we define $F(\mathcal{R}[t, t^{-1}])$.

In defining $F(\mathcal{R}[t, t^{-1}])$ we first define k x k matrices $T_i$ ($i = 2,\ldots,k$). $T_i$ has t in the first i-1 diagonal entries, 1 in the remaining diagonal entries, 1-t in the $i^{th}$ row prior to the diagonal, and zeros everywhere else. Notice that $T_2 = T$ when k is 2. It is easy to see that $gp\langle T_2, \ldots, T_k \rangle$ is free abelian of rank k-1, but all that is necessary is that it be abelian (in generalizing Lemma 5 for arbitrary ranks). We can now define $F(\mathcal{R}[t, t^{-1}])$ as

$$F(\mathcal{R}[t, t^{-1}]) = gp\langle M_1, M_2 T_2, M_3 T_3, \ldots, M_k T_k \rangle.$$

In dealing with the Burnside problem for k generators we can now follow the procedure for two generator groups. The proof that $F(\mathcal{R})$ is free metabelian of rank k and that $F(S)$ is the free metabelian Burnside group of exponent q and rank k should follow as previously. As before, we have the commutative square where the horizontal maps $\alpha$ and $\underline{\alpha}$ come from the ring homomorphism $\mathcal{R} \to S$ and the vertical maps are those determined by sending t to the identity.

$$\begin{array}{ccc} F(\mathcal{R}[t, t^{-1}]) & \xrightarrow{\alpha} & F(S[t, t^{-1}]) \\ \downarrow & & \downarrow \\ F(\mathcal{R}) & \xrightarrow{\underline{\alpha}} & F(S) \end{array}$$

When $k = 2$, our notation reduces to the same as in prior sections. $\mathcal{R} = \mathcal{R}(2) = Z[x, x^{-1}, y, y^{-1}]$ and $F(\mathcal{R})$, $F(S)$, and $F(\mathcal{R}[t, t^{-1}])$ are as before. $F(\mathcal{R}) = gp\langle M_1, M_2 \rangle$ and $F(\mathcal{R}[t,t^{-1}]) = gp\langle M_1, M_2 T \rangle$.

Details for k generators need to be worked out.

## Appendix

Lemma 1(ii): Each element M of $F(\mathcal{R})$ has the form $M = uI + N$, where I is the identity matrix, u is a positive unit in $\mathcal{R}$, and $N = [\lambda_i v]$ is an r-square matrix whose $i^{th}$ row is $\lambda_i v$. The $\lambda_i$ in $\mathcal{R}$ satisfy $\lambda_1(1-x_1) + \cdots + \lambda_r(1-x_r) = 1- u$.

Proof: It is a simple matter to check that $vM = v$ for all M in $F(\mathcal{R})$. (i.e. Verify this for each generator in $F(\mathcal{R})$.) Also verify that each generator and its inverse satisfies Lemma 1(ii). Having verified Lemma 1(ii) for words of length 1, we use induction on the length of a word in $F(\mathcal{R})$.

Suppose $W_1W_2$ is an element of $F(\mathcal{R})$, where $W_1 = u_1 + [\lambda_i v]$, $W_2 = u_2 + [\delta_i v]$ and $\lambda_1(1-x_1) + \cdots + \lambda_r(1-x_r) = 1 - u_1$, $\delta_1(1-x_1) + \cdots + \delta_r(1-x_r) = 1 - u_2$. Then using the fact that $vM = v$ for all $M$ in $F(\mathcal{R})$,

$W_1W_2 = (u_1 + [\lambda_i v])(u_2 + [\delta_i v]) = u_1(u_2 + [\delta_i v]) + [\lambda_i v] = u_1u_2 + u_1[\delta_i v]) + [\lambda_i v]$. But $(u_1\delta_1 + \lambda_1)(1-x_1) + \cdots + (u_1\delta_r + \lambda_r)(1-x_r) = u_1(1 - u_2) + (1 - u_1) = 1 - u_1u_2$. Thus the result holds for $W_1W_2$

301 El Monte Dr
Santa Barbara Ca. 93109
email: sb7@me.com

e4e4                                22